\documentclass[11pt]{amsart}
\usepackage{amsmath, amssymb, amscd, mathrsfs, url, pinlabel,verbatim}
\usepackage[pagebackref,hidelinks]{hyperref}
\usepackage[margin=1.25in,marginparwidth=1in,centering,letterpaper,dvips]{geometry}
\usepackage{color,dcpic,latexsym,graphicx,epstopdf,comment}
\usepackage[all]{xy}
\usepackage[dvipsnames]{xcolor}
\usepackage{tikz,tikz-cd,pgfplots}
\usepackage{upquote,tabularx,textcomp, booktabs}
\usepackage[shortlabels]{enumitem}
\usepackage[color=blue!20!white,textsize=tiny]{todonotes}

\title{3-manifold decision problems and Riemann hypotheses}

\author{John A. Baldwin}
\address{Department of Mathematics \\ Boston College}
\email{john.baldwin@bc.edu}

\author[Steven Sivek]{Steven Sivek}
\address{Department of Mathematics\\Imperial College London}
\email{s.sivek@imperial.ac.uk}

\makeatletter
\newtheorem*{rep@theorem}{\rep@title}
\newcommand{\newreptheorem}[2]{%
\newenvironment{rep#1}[1]{%
 \def\rep@title{#2 \ref{##1}}%
 \begin{rep@theorem}}%
 {\end{rep@theorem}}}
\makeatother

\newtheorem {theorem}{Theorem}
\newreptheorem{theorem}{Theorem}
\newtheorem {lemma}[theorem]{Lemma}
\newtheorem {proposition}[theorem]{Proposition}
\newtheorem {corollary}[theorem]{Corollary}

\newtheorem {hypothesis}[theorem]{Hypothesis}

\numberwithin{equation}{section}
\numberwithin{theorem}{section}

\theoremstyle{definition}

\newtheorem{remark}[theorem]{Remark}
\newtheorem*{remark*}{Remark}

\newlist{pcases}{enumerate}{1}
\setlist[pcases]{
  label=\bf{Case~\arabic*:}\protect\thiscase.~,
  ref=\arabic*,
  align=left,
  labelsep=0pt,
  leftmargin=0pt,
  labelwidth=0pt,
  parsep=0pt
}
\newcommand{\case}[1][]{%
  \if\relax\detokenize{#1}\relax
    \def\thiscase{}%
  \else
    \def\thiscase{~#1}%
  \fi
  \item
}

\newcommand{\Z}{\mathbb{Z}}

\newcommand{\C}{\mathbb{C}}

\newcommand{\F}{\mathbb{F}}
\newcommand{\Q}{\mathbb{Q}}

\newcommand{\Hom}{\operatorname{Hom}}

\newcommand{\cO}{\mathcal{O}}
\newcommand{\fp}{\mathfrak{p}}
\newcommand{\fq}{\mathfrak{q}}

\newcommand\SL{\mathrm{SL}}

\DeclareMathOperator{\poly}{poly}

\DeclareMathOperator{\frob}{Frob}
\DeclareMathOperator{\Gal}{Gal}
\DeclareMathOperator{\Li}{Li}

\usetikzlibrary{decorations.markings,arrows,calc,intersections,patterns}
\tikzset{every picture/.style=thick}
\tikzset{link/.style = { white, double = black, line width = 1.75pt, double distance = 1.25pt, looseness=1.75 }}
\tikzset{crossing/.style = {draw, circle, dotted, minimum size=0.5cm, inner sep=0, outer sep=0}}
\pgfplotsset{compat=1.12}

\begin{document}

\begin{abstract}
Beginning with Kuperberg's work on unknot recognition \cite{kuperberg}, several fundamental decision problems in 3-manifold topology have been shown to be in coNP assuming the generalized Riemann hypothesis. In this note, we explain how recent developments around the quasi-Riemann hypothesis yield unconditional versions of these results, leading for example to the first unconditional proof that 3-sphere recognition is in coNP. 
\end{abstract}

\maketitle

\section{Introduction}
A decision problem lies in the complexity class NP if given an input for which the answer is \emph{yes}, there exists a certificate of this fact whose length and verification time are polynomial in the size of the input; it lies in coNP if instead there exists such a certificate when the answer is \emph{no}. Several fundamental questions in 3-manifold topology are known to lie in NP. These include unknot recognition \cite{hass-lagarias-pippenger}, torus knot recognition \cite{bs-complexity}, satellite knot recognition \cite{haraway-hoffman}, recognition of any fixed knot \cite{lackenby-reidemeister}, and 3-sphere recognition \cite{ivanov, schleimer}, among others.

In \cite{kuperberg}, Kuperberg gave a beautiful proof that unknot recognition is also in coNP, assuming the generalized Riemann hypothesis (GRH). His certificate of knottedness relies on a deep gauge-theoretic result of Kronheimer and Mrowka \cite{km-su2}, which asserts that if a knot $K\subset S^3$ is not the unknot, then there exists a representation \[ \pi_1(S^3\setminus K) \to \SL(2,\C) \] with nonabelian image.  Using work of Koiran \cite{koiran}, he deduced that there is a nonabelian representation \[ \rho: \pi_1(S^3\setminus K) \to \SL(2,\F_p) \] for some prime $p$, and that by work of Lagarias--Odlyzko \cite{lagarias-odlyzko}, GRH allows us to take $p$ with at most polynomially many digits in the number of crossings in any planar diagram for $K$ (or the number of tetrahedra in any triangulation of $S^3\setminus K$). The certificate of knottedness then consists of matrices $\rho(g_i)$, where $g_i$ ranges over generators in a presentation of $\pi_1(S^3\setminus K)$ coming from the input.

The same ideas have been used in combination with additional results from gauge theory \cite{km-excision, zentner} on a variety of decision problems, which we summarize here.

\begin{theorem} \label{thm:coNP-conditional}
The following problems are in coNP, assuming GRH:
\begin{enumerate}
\item unknot recognition \cite{kuperberg}; \label{i:unknot}
\item 3-sphere recognition \cite{zentner}; \label{i:3-sphere}
\item torus knot recognition and hyperbolic knot recognition \cite{bs-complexity}; \label{i:torus-hyperbolic}
\item hyperbolicity of cusped 3-manifolds \cite{haraway-hoffman}. \label{i:hyperbolicity}
\end{enumerate}
The input can be either a planar knot diagram in cases~\eqref{i:unknot} and \eqref{i:torus-hyperbolic}, or a triangulation in cases~\eqref{i:unknot}, \eqref{i:3-sphere}, and \eqref{i:hyperbolicity}.
\end{theorem}

In each case, GRH is used to ensure that the certificate has polynomial length in the input, but once we know that a certificate of this length exists, the fact that it can be verified in polynomial time is unconditionally true.  (Haraway and Hoffman did not state case~\eqref{i:hyperbolicity} in terms of GRH, but they proved that hyperbolicity is in coNP if 3-sphere recognition is.)  Several of these problems have since been proved to be in coNP unconditionally by different means, as detailed in Remark~\ref{rmk:unconditional} below, but the proofs relying on GRH are still interesting because the corresponding certificates are arguably much simpler.

While writing our paper \cite{bs-complexity}, we made the unfortunate decision to treat the number-theoretic input as a black box, taking ``GRH'' as something mysterious but mostly harmless.  This is fairly common practice (many of us use ``ZFC'' in the same way) but not exactly a sign of understanding on our part---for example, it would have been good to know if partial progress toward GRH might actually make the above theorems unconditional, as explicitly suggested in the introduction to \cite{kuperberg}.  The purpose of this note is to atone for our sins by seeing what exactly we might need from GRH to prove these theorems unconditionally.  We do not claim any novelty in our arguments, but we have spelled them out explicitly so that topologists might see how this circle of ideas can be applied.

Our starting point is the following explicit claim, called the \emph{quasi-Riemann hypothesis} for Dedekind zeta functions.  (The earliest reference we could find for this name is \cite{turan}, but the conjecture itself goes back at least to Hardy and Littlewood \cite[Hypothesis~R]{hardy-littlewood}.)  By way of notation, we write $\cO_K$ for the ring of integers in a number field $K$, and $N_{K/\Q}(I) = [\cO_K : I]$ denotes the norm of an ideal $I \subset \cO_K$.

\begin{hypothesis} \label{hyp:weak-GRH}
There is a universal constant $\epsilon > 0$ such that for every number field $K$, the Dedekind zeta function
\[ \zeta_K(s) = \sum_{I \subset \cO_K} \frac{1}{N_{K/\Q}(I)^s} = \prod_{\mathfrak{p} \subset \cO_K \text{ prime}} \left(1 - \frac{1}{N_{K/\Q}(\mathfrak{p})^s}\right)^{-1} \]
has no zeroes with $\operatorname{Re}(s) > 1-\epsilon$.
\end{hypothesis}

The purpose of this note is to explain how the arguments of \cite{kuperberg} can also be used to prove the following.

\begin{theorem} \label{thm:weak-GRH-implies-coNP}
Each of the decision problems in Theorem~\ref{thm:coNP-conditional} is in coNP, modulo Hypothesis~\ref{hyp:weak-GRH}.
\end{theorem}

The results of Theorem~\ref{thm:coNP-conditional} all rely on the existence of a suitably small prime $p$, which Koiran \cite{koiran} finds (assuming GRH) by appealing to Lagarias and Odlyzko's \emph{effective} version of the Chebotarev density theorem \cite{lagarias-odlyzko}.  Here the term ``effective'' means that it describes not only the asymptotic growth of some quantity but also a bound on the error term.  Frei, Sofos, and Widmer \cite{frei-sofos-widmer} have recently shown that Hypothesis~\ref{hyp:weak-GRH} implies an effective version whose error bounds are weaker than the ones implied by GRH (see \S\ref{sec:chebotarev} for details), but these turn out to be sharp enough to prove the existence of our sought-after prime $p$.  This existence is established in Theorem~\ref{thm:goal}, which can be used verbatim in place of \cite[Theorem~3.2]{kuperberg} to prove Theorem~\ref{thm:weak-GRH-implies-coNP}. 

Cureton and Thorner \cite{cureton-thorner} have recently communicated a Lean-formalized proof of Hypothesis~\ref{hyp:weak-GRH} with $\epsilon = \frac{1}{32}$, building on the formalized proof announced by OpenAI \cite{oai1,oai2} of the quasi-Riemann hypothesis for Dirichlet $L$-functions.  This implies via Theorem~\ref{thm:weak-GRH-implies-coNP} that the results of Theorem~\ref{thm:coNP-conditional} no longer rely on any open conjectures: the named decision problems are unconditionally in coNP.  However, at the time of writing it is not clear to us that any human has understood the proof given in \cite{cureton-thorner}, so following the lead of \cite{frei-sofos-widmer} we will continue to describe these results as being conditional on Hypothesis~\ref{hyp:weak-GRH}.

\begin{remark}\label{rmk:unconditional}
Table~\ref{tab:complexity} summarizes what is known about each of the recognition problems of Theorem~\ref{thm:coNP-conditional}.  In the NP column, Baroni and Lackenby have announced forthcoming work proving that hyperbolic knot recognition is in NP.  For the unconditionally coNP column, some of these were already known by different methods: this includes unknot recognition, as proved by Lackenby \cite{lackenby-knottedness}, and also by Agol in unpublished work from around 2002; and hyperbolic knot recognition, by Haraway and Hoffman \cite[Theorem~2]{haraway-hoffman}.  As far as we are aware, however, there are no other known strategies in the literature for 3-sphere recognition or hyperbolicity of cusped 3-manifolds.
\end{remark}

\begin{table}
\label{tab:complexity}
\begin{tabular}{lccc}
\toprule
Problem
& NP
& coNP mod GRH
& coNP \\
\midrule
Unknot recognition
 & \cite{hass-lagarias-pippenger}
 & \cite{kuperberg}
 & \cite{lackenby-knottedness} \\
\addlinespace
Torus knot recognition
& \cite{bs-complexity}
& \cite{bs-complexity}
& Theorem~\ref{thm:weak-GRH-implies-coNP} \\
\addlinespace
Hyperbolic knot recognition
& (forthcoming)
& \cite{bs-complexity}
& \cite{haraway-hoffman} \\
\addlinespace
Hyperbolicity of cusped $3$-manifolds
& ?
& \cite{zentner}$+$\cite{haraway-hoffman}
& Theorem~\ref{thm:weak-GRH-implies-coNP} \\
\addlinespace
$3$-sphere recognition
 & \cite{ivanov,schleimer}
 & \cite{zentner}
 & Theorem~\ref{thm:weak-GRH-implies-coNP} \\
\bottomrule
\end{tabular}
\vspace{8pt}
\caption{References for the complexity of various 3-manifold decision problems.  The results attributed to Theorem~\ref{thm:weak-GRH-implies-coNP} assume only Hypothesis~\ref{hyp:weak-GRH}, a proof of which was claimed in \cite{cureton-thorner}.}
\end{table}

\begin{remark}
Hypothesis~\ref{hyp:weak-GRH} also implies nontrivial bounds on the least size of nonabelian finite quotients of many knot and 3-manifold groups, exactly as in \cite[Theorem~5.1]{kuperberg}.  See Theorem~\ref{thm:finite-cover} for a precise statement.
\end{remark}

\subsection*{AI statement}

No generative AI was used at any point in the writing of this note.

\subsection*{IA statement}

As we were writing this note, Ian Agol informed us that he was also aware that Hypothesis~\ref{hyp:weak-GRH} suffices to prove an unconditional version of Theorem~\ref{thm:coNP-conditional}.

\subsection*{Acknowledgments} JAB was supported by NSF Grant DMS-2506250 and a Simons Fellowship.  SS was supported by the Engineering and Physical Sciences Research Council [grant number UKRI1016]. No data were created or analyzed in this work.

\section{Background} \label{sec:background}

The results of Theorem~\ref{thm:coNP-conditional} all rely on the same technique.  Given a decision problem, we construct a group $G$ with a modestly-sized presentation---this could be a knot group $\pi_1(S^3 \setminus K)$, or a 3-manifold group $\pi_1(Y)$---and then the representation variety
\[ \Hom(G, \SL(2,\C)) \]
is a complex affine algebraic set that can be described by a controlled number of variables and relations.  If we want to find a representation $\rho$ with non-abelian image, this can be achieved at minimal extra cost by the Rabinowitsch trick: given polynomials $w_1,\dots,w_k$ in the input variables, we can impose the condition $(w_1,\dots,w_k) \neq (0,\dots,0)$ by adding $k$ extra variables $t_i$ and the equation $\sum_i t_iw_i = 1$.  Applying this trick to the entries of $\rho(g_ig_j)-\rho(g_jg_i)$ as $g_i,g_j$ range over pairs of generators of $G$, we thus get a system of polynomial equations with integral coefficients
\[ f_i(x_1,\dots,x_n) = 0, \qquad 1 \leq i \leq m, \]
with the properties that
\begin{enumerate}
\item the integers $m,n$ and the coefficients of the $f_i$ are all bounded in terms of the input;
\item and the system has a solution if and only if the desired representation $\rho$ exists.
\end{enumerate}

If this system has a solution, then Kuperberg \cite[Theorem~3.1]{kuperberg} asserts that by work of Koiran \cite[Theorems~4 and 7]{koiran}, there is an explicit solution in terms of a root of a single irreducible polynomial
\[ f \in \Z[x], \]
whose degree and coefficients are suitably bounded by $m$, $n$, and the degrees and coefficients of the $f_i$.  In particular, if $f(x)=0$ has a solution modulo $p$ then so will our original system of equations, so our goal will be to find such a $p$ that is not too large.  In other words, we just need to prove the following strengthening of \cite[Theorem~3.2]{kuperberg}.

\begin{theorem} \label{thm:goal}
Assume Hypothesis~\ref{hyp:weak-GRH}.  Let $f \in \Z[x]$ be an irreducible polynomial of degree $d \geq 1$, and let $||f||$ denote the maximum absolute value among the coefficients of $f$.  Then there is a prime
\[ p = \poly(d, \log ||f||) \]
such that $f$ has a root modulo $p$.
\end{theorem}

We note that Theorem~\ref{thm:goal} will then imply the following corollaries, by exactly the same arguments as for the corresponding theorems in \cite{kuperberg}.

\begin{theorem}[{cf.~\cite[Theorem~3.4]{kuperberg}}] \label{thm:affine-group}
Assume Hypothesis~\ref{hyp:weak-GRH}, and let $G$ be an affine algebraic group over $\Z$.  Let $\Gamma$ be a discrete group with a finite presentation of length $\ell$.  If there is a homomorphism $\Gamma \to G(\C)$ with nonabelian image, then there is also a homomorphism $\Gamma \to G(\F_p)$ with nonabelian image for some prime $p$ with $\log(p) = \poly(\ell)$.
\end{theorem}

\begin{theorem}[{cf.~\cite[Theorem~5.1]{kuperberg}}] \label{thm:finite-cover}
Assume Hypothesis~\ref{hyp:weak-GRH}.  Let $Y$ be any of the following:
\begin{enumerate}
\item a finite-volume hyperbolic 3-manifold,
\item the complement of a nontrivial knot $K \subset S^3$,
\item or an integer homology 3-sphere not homeomorphic to $S^3$.
\end{enumerate}
Suppose that $Y$ admits a triangulation with $n$ tetrahedra.  Then there is a regular nonabelian cover $\tilde{Y} \to Y$ of degree at most $e^{\poly(n)}$.
\end{theorem}

\begin{proof}
There is a homomorphism $\pi_1(Y) \to \SL(2,\C)$ with nonabelian image, by \cite[Proposition~3.1.1]{culler-shalen}, \cite{km-su2}, or \cite{zentner} respectively.  Since $\pi_1(Y)$ has a presentation of length $\poly(n)$, Theorem~\ref{thm:affine-group} gives a nonabelian
\[ \rho : \pi_1(Y) \to \SL(2,\F_p) \]
for some prime $p$ with $\log(p) = \poly(n)$.  Then $\pi_1(Y)$ surjects onto the image of $\rho$, which is a nonabelian group of order at most
\[ \left| \SL(2,\F_p) \right| = p(p^2-1) < p^3, \]
and the regular cover of $Y$ corresponding to $\ker(\rho)$ has degree $|\operatorname{im}(\rho)|$, which is less than $p^3 = (e^{\poly(n)})^3 = e^{\poly(n)}$.
\end{proof}

Returning to Theorem~\ref{thm:goal}, the original version stated in \cite{kuperberg}, which is attributed to Lagarias--Odlyzko \cite{lagarias-odlyzko} and Weinberger \cite{weinberger}, assumed GRH.  Our proof follows the same lines, with GRH replaced by Hypothesis~\ref{hyp:weak-GRH}, but we spell out the details much more explicitly here.

The strategy is as follows.  It is hard to find an explicit choice of $p$ that works.  However, many classical results in number theory describe the \emph{statistical} behavior of prime numbers, such as the prime number theorem: if we define
\[ \Li(x) = \int_2^x \frac{dt}{\log t}, \]
then the number of primes less than a given $x$ is asymptotic to $\Li(x)$.  This guarantees that prime numbers exist, but it doesn't tell us how far we must search before we can be certain that we will find one.  (This is strictly for motivation; we do actually know that $2,3,5,7,\dots$ are prime.)  In order to get effective bounds, we need to estimate the size of the error term, and the Riemann hypothesis tells us that this error should be $O(x^{1/2+\epsilon})$ for any $\epsilon > 0$.

In order to prove Theorem~\ref{thm:goal}, then, we define a suitable counting function: given a prime $p$, we let
\[ W(p) = \#\{ a \in \F_p \mid f(a) = 0 \}, \]
and we want to know that $W(p)$ is positive for some not-too-large value of $p$.  Equivalently, we want to find some relatively small $x$ such that the sum
\[ \sum_{p \leq x} W(p) \]
is strictly positive.  If we count roots of $f$ mod $p$ with multiplicity, then we can think of $W(p)$ as actually counting the number of linear factors in the prime factorization of $f$, where $f$ is viewed as an element of $\F_p[x]$.  (We mostly don't have to worry about multiplicity: if $p$ does not divide the discriminant $\Delta_f$, then $f$ has no repeated factors mod $p$.)  We are thus led to the following classical theorem of Dedekind and Kummer.

\begin{theorem}[{\cite[Proposition I.8.3]{neukirch}}] \label{thm:dedekind-kummer}
Let $f \in \Z[x]$ be monic and irreducible.  Let $\theta$ be a root of $f$, and let $K = \Q(\theta) \cong \Q[x]/\langle f\rangle$.  If $p$ does not divide the index $[\cO_K  : \Z[\theta]]$, and the mod $p$ reduction $\bar{f} \in \F_p[x]$ of $f$ has irreducible factorization
\[ \bar{f} = \bar{f}_1^{e_1} \dots \bar{f}_k^{e_k}, \qquad \bar{f}_i \in \F_p[x], \]
then the ideal $\langle p \rangle \subset \cO_K$ has prime factorization
\[ p\cO_K = \prod_{i=1}^k \fq_i^{e_i}, \qquad \fq_i = \langle p, f_i(\theta)\rangle \]
where each $f_i \in \Z[x]$ is a lift of $\bar{f}_i \in \F_p[x]$.  Each prime $\fq_i$ has norm $N_{K/\Q}(\fq_i) = p^{\deg \bar{f}_i}$.
\end{theorem}

See also \cite[Theorem~8.2]{stevenhagen} for a slightly more general version.

Theorem~\ref{thm:dedekind-kummer} tells us that for a given monic, irreducible $f$, if we set $K = \Q[x]/\langle f(x)\rangle$, then the sum $\sum_{p \leq x} W(p)$ is closely related to the count
\begin{equation} \label{eq:count-prime-ideals}
\# \{ \fq \subset \cO_K \mid N_{K/\Q}(\fq) \leq x \}
\end{equation}
of prime ideals of norm not exceeding $x$, and more precisely to the count of those ideals whose norm is a rational prime.  We must therefore try to estimate the count \eqref{eq:count-prime-ideals}, and then also bound the smaller count of prime ideals whose norms are prime powers $p^e$ for some $e \geq 2$.

\section{The effective Chebotarev density theorem} \label{sec:chebotarev}

In this section we explain the input needed from \cite{frei-sofos-widmer} to prove Theorem~\ref{thm:goal}.  Strictly speaking, we will only need Corollary~\ref{cor:landau}, which is easier to state than the prerequisites that we are about to explain.  We encourage the topologists reading this to stick with us and learn about the Chebotarev density theorem anyway, rather than skipping ahead.

Let $K$ be a number field, and let $L$ be a normal extension of $K$ with Galois group
\[ G = \Gal(L/K). \]
A prime ideal $\fp \subset \cO_K$ is said to be \emph{unramified} in $L$ if the ideal $\fp \cO_L \subset \cO_L$ has prime factorization
\[ \fp \cO_L = \fq_1 \fq_2 \cdot \ldots \cdot \fq_k \]
without any repeated factors.  Then for each ideal $\fq = \fq_i$, there is a uniquely defined Frobenius automorphism
\[ \frob_\fq \in \Gal(L/K) \]
fixing $\fq$ setwise and satisfying $\frob_\fq(t) \equiv t^{\#|\cO_K/\fp|} \pmod{\fq}$ for all $t \in \cO_L$.  In general $\frob_\fq$ depends on our choice of prime $\fq$ lying over $\fp$, but only up to conjugacy.  We use the Artin symbol $\left[\frac{L/K}{\fp}\right]$ to denote this conjugacy class in $\Gal(L/K)$.

Following the notation of Lagarias and Odlyzko \cite{lagarias-odlyzko}, we set
\[ \pi_C(x, L/K) = \# \left\{ \fp \subset \cO_K \,\middle\vert\, \begin{array}{c} \fp \text{ unramified in } L \\[0.25em] \left[\frac{L/K}{\fp}\right] = C \\[0.25em] N_{K/\Q}(\fp) \leq x \end{array} \right\}. \]
The Chebotarev density theorem says that for any conjugacy class $C$ in $G = \Gal(L/K)$, we have the asymptotic relation
\[ \pi_C(x,L/K) \sim \frac{\#C}{\#G} \operatorname{Li}(x) \]
as $x \to \infty$.  The classical version of this theorem says nothing about the error terms, but Lagarias and Odlyzko's \emph{effective} Chebotarev density theorem \cite[Theorem~1.1]{lagarias-odlyzko} bounds the error from above under the assumption of GRH.  Frei, Sofos, and Widmer \cite{frei-sofos-widmer} have shown that the proof applies equally well under the assumption of Hypothesis~\ref{hyp:weak-GRH}:

\begin{theorem}[{\cite[Theorem~1.1]{frei-sofos-widmer}}] \label{thm:effective-chebotarev}
Assume Hypothesis~\ref{hyp:weak-GRH}.  There exists a universal constant $c > 0$ such that for every $x > 2$, if $L$ is a normal extension of the number field $K$ with Galois group $G = \Gal(L/K)$, then we have
\[ \left| \pi_C(x, L/K) - \frac{\#C}{\#G} \operatorname{Li}(x) \right| \leq c \cdot \frac{\#C}{\#G} x^{1-\epsilon}\log (\Delta_L^{1/\epsilon} x^{[L:\Q]}). \]
\end{theorem}

\begin{remark}
The theorem stated in \cite{lagarias-odlyzko} corresponds to the case $\epsilon = \frac{1}{2}$.
\end{remark}

Theorem~\ref{thm:effective-chebotarev} quickly implies an effective version of Landau's prime ideal theorem, which is what we will actually need to prove Theorem~\ref{thm:goal}.

\begin{corollary} \label{cor:landau}
Assume Hypothesis~\ref{hyp:weak-GRH}.  Then there exists a universal constant $c > 0$ with the following property: if $f \in \Z[x]$ is an irreducible monic polynomial of degree $d \geq 1$, and $K = \Q[x] / \langle f(x)\rangle$, then
\[ \left| \#\left\{ \fq \subset \cO_K \,\middle\vert\, N_{K/\Q}(\fq) \leq x \right\} - \Li(x) \right| \leq c \cdot x^{1-\epsilon} \log(\Delta_K^{1/\epsilon}x^d) \]
for every $x > 2$.
\end{corollary}

\begin{proof}
We take $L = K$, which is trivially a normal extension of $K$.  Then $\Gal(L/K)$ is the trivial group, and we apply Theorem~\ref{thm:effective-chebotarev} to its unique conjugacy class.
\end{proof}

Corollary~\ref{cor:landau} says that the quantity \eqref{eq:count-prime-ideals}, which approximates the total number of roots our monic polynomial $f$ has modulo the various primes $p \leq x$, is at least
\[ \Li(x) - c \cdot x^{1-\epsilon}\log(\Delta_K^{1/\epsilon}x^d). \]
Thus we have two remaining goals, at least in the proof of Theorem~\ref{thm:goal} for monic $f$: one is to determine how many of the prime ideals in $\cO_K$ actually correspond to roots mod $p$, and then the other is to the bound error term $x^{1-\epsilon}\log(\Delta_K^{1/\epsilon}x^d)$.

In order to understand the discriminant $\Delta_K$ appearing in this error term, we apply the identity
\begin{equation} \label{eq:discriminant-conductor}
\Delta_f = [\cO_K : \Z[\theta]]^2 \cdot \Delta_K,
\end{equation}
see e.g.~\cite[Equation~(7.7)]{stevenhagen}.  Here the term $\Delta_f$ denotes the absolute value of the product
\[ \prod_{1 \leq i < j \leq d} (\theta_i - \theta_j)^2, \]
where the $\theta_i$ are the roots of $f$ in $\C$.

\begin{lemma}[cf.\ {\cite[Equation~(3)]{kuperberg}}] \label{lem:discriminant-bound}
Let $f \in \Z[x]$ be a monic polynomial of degree $d \geq 1$.  The discriminant $\Delta_f$ satisfies
\[ \log \Delta_f < 2d\big( \log(d+1) + \log ||f|| \big). \]
\end{lemma}

\begin{proof}
We can evaluate $\Delta_f$ as the resultant of $f$ and $f'$, which is defined as the determinant of a certain $(2d-1) \times (2d-1)$ matrix $M$, the Sylvester matrix of $f$ and $f'$.  The first $d-1$ rows of $M$ each contain shifts of the coefficient vector of $f$: there are at most $d+1$ nonzero entries, each bounded in absolute value by $||f||$, so as row vectors these have length at most $||f||\sqrt{d+1}$.  Likewise, the remaining $d$ rows are all shifts of the coefficient vector of $f'$, which has at most $d$ nonzero entries bounded by $d||f||$, so as row vectors these have length at most $d||f||\cdot\sqrt{d}$.  Then $\Delta_f = \left|\det(M)\right|$ is at most the product of the lengths of the rows of $M$:
\begin{align*}
\Delta_f &\leq \left(||f||\sqrt{d+1}\right)^{d-1} \cdot \left(||f||d^{3/2}\right)^d \\
& < ||f||^{2d-1} (d+1)^{2d-1/2}.
\end{align*}
The lemma follows immediately.
\end{proof}

\section{Counting roots mod $p$ of monic polynomials} \label{sec:monic}

Let $f \in \Z[x]$ be monic and irreducible, of degree $d \geq 1$, and with some root $\theta \in \C$.  Letting $K = \Q(\theta)$ as in Theorem~\ref{thm:dedekind-kummer}, we now attempt to classify the primes $\fq \subset \cO_K$ of norm at most $x$.  Each such $\fq$ lies over some rational prime $p$, meaning that $\fq \cap \Z = p\Z$, and then it is one of the factors in the prime factorization
\[ p\cO_K = \fq_1 \dots \fq_k. \]
Since the relative norm of an ideal is multiplicative, we have
\begin{equation} \label{eq:factor-norm-multiplicative}
\prod_{i=1}^k N_{K/\Q}(\fq_i) = N_{K/\Q}(p\cO_K) = p^d,
\end{equation}
so each $N_{K/\Q}(\fq_i)$ is a positive power of $p$; there are at most $d$ factors $\fq_i$; and
\[ p \leq N_{K/\Q}(\fq) \leq x. \]
Thus each prime ideal $\fq \subset \cO_K$ of norm at most $x$ lies over a prime $p \leq x$ satisfying exactly one of the following:
\begin{enumerate}
\item $p$ does not divide $\Delta_f$, and $N_{K/\Q}(\fq) = p$; \label{i:root-mod-p}
\item $p$ does not divide $\Delta_f$, and $N_{K/\Q}(\fq) = p^e$ for some $e \geq 2$; \label{i:not-root-mod-p}
\item $p$ divides $\Delta_f$. \label{i:ramified}
\end{enumerate}
We consider each of these cases separately, in order to bound case~\eqref{i:root-mod-p} from below as done in the proof of \cite[Corollary~1]{koiran}.

In case~\eqref{i:root-mod-p}, we recall from \eqref{eq:discriminant-conductor} that $p$ cannot divide $[\cO_K : \Z[\theta]]$, so Theorem~\ref{thm:dedekind-kummer} says that such ideals $\fq$ are in bijection with the roots of $f$ mod $p$.

In case~\eqref{i:not-root-mod-p} we have $p^2 \leq N_{K/\Q}(\fq) \leq x$, so $p \leq \sqrt{x}$.  For any fixed $p$, the relation \eqref{eq:factor-norm-multiplicative} says that there are at most $d/2$ such primes $\fq$ over $p$, so overall this case contributes at most $d\sqrt{x}/2$ prime ideals to the total.

Finally, in case~\eqref{i:ramified}, there are at most $\log_2 \Delta_f < 2\log \Delta_f$ primes $p$ dividing $\Delta_f$, and for each of these the ideal $p\cO_K$ contributes at most $d$ prime factors $\fq$ of norm at most $x$, so in total this case contributes at most $2d\cdot \log \Delta_f$ prime ideals.

Isolating case~\eqref{i:root-mod-p} from the rest, we have shown that
\begin{equation} \label{eq:lower-bound-unconditional}
\sum_{\substack{p \leq x \\ p \nmid \Delta_f}} W(p) \geq \#\{\fq \subset \cO_K \mid N_{K/\Q}(\fq) \leq x\}  - \tfrac{1}{2}dx^{1/2} - 2d\log \Delta_f.
\end{equation}

\begin{proposition} \label{prop:goal-monic}
Assuming Hypothesis~\ref{hyp:weak-GRH}, there is a universal constant $c > 0$ with the following property.  Let $f \in \Z[x]$ be a monic, irreducible polynomial of degree $d \geq 1$.  Then
\[ \sum_{\substack{p \leq x \\ p \nmid \Delta_f}} W(p) \geq \Li(x) - c x^{1-\epsilon}\log(\Delta_f^{1/\epsilon}x^d) - d\left(\tfrac{1}{2}x^{1/2}+2\log\Delta_f\right) \]
for all $x > 2$.
\end{proposition}

\begin{proof}
Let $c$ be the constant provided by the effective prime ideal theorem, Corollary~\ref{cor:landau}, so that
\[ \#\{\fq \subset \cO_K \mid N_{K/\Q}(\fq) \leq x\} \geq \Li(x) - c\cdot x^{1-\epsilon} \log(\Delta_K^{1/\epsilon} x^d). \]
We can replace $\Delta_K$ with $\Delta_f$, which is at least as large by \eqref{eq:discriminant-conductor}, and combine this with \eqref{eq:lower-bound-unconditional} to conclude.
\end{proof}

\section{Proof of Theorem~\ref{thm:goal}} \label{sec:proof}

In our work towards proving Theorem~\ref{thm:goal}, we have assumed until now that the polynomial $f \in \Z[x]$ is monic, but the general case quickly reduces to this one.  Suppose that
\[ f(x) = a_d x^d + a_{d-1} x^{d-1} + \dots + a_0, \]
where $a_d \neq 0$.  We consider the related polynomial
\[ g(x) = a_d^{d-1} f\left(\frac{x}{a_d}\right) = x^d + \sum_{i=0}^{d-1} a_d^{d-1-i} a_i x^i, \]
which is monic and which has integer coefficients as well.  Then $\theta$ is a root of $f(x)$ if and only if $a_d\theta$ is a root of $g(x)$, and this is true not just over $\C$ but also over any field $\F_p$ as long as $a_d$ is invertible mod $p$.  Thus if we want to find a prime $p$ such that $f$ has a root mod $p$, it suffices to find $p \nmid a_d$ such that $g$ has a root mod $p$.

We now observe that $\deg g = \deg f = d$, and also that 
\[ ||g|| \leq |a_d|^{d-1} ||f|| \leq ||f||^d \qquad \Longrightarrow \qquad \log ||g|| \leq d \cdot \log ||f||. \]
Since $|a_d| \leq ||f||$, the number of primes dividing $a_d$ is at most
\[ \log_2 |a_d| \leq 2 \log |a_d| \leq 2\log ||f||, \]
and each of them has at most $d$ prime ideal factors in $\cO_K$ where $K = \Q[x]/\langle g\rangle$.  Thus in this case we can weaken the bound of Proposition~\ref{prop:goal-monic} slightly, by removing all the prime ideals that divide some rational prime divisor of $a_d$, to get
\[ \sum_{\substack{p \leq x \\ p \nmid \Delta_g \\ p \nmid a_d}} W(p) \geq \Li(x) - c x^{1-\epsilon}\log(\Delta_g^{1/\epsilon}x^d) - d\left(\tfrac{1}{2}x^{1/2}+2\log\Delta_g + 2\log ||f||\right). \]
Applying Lemma~\ref{lem:discriminant-bound} to the monic polynomial $g$ gives us
\begin{equation} \label{eq:delta-g}
\begin{aligned}
\log \Delta_g &< 2d\big( \log(d+1) + \log ||g|| \big) \\
&\leq 2d\big( \log(d+1) + d\log ||f|| \big),
\end{aligned}
\end{equation}
so collectively we have shown the following.

\begin{proposition} \label{prop:goal-general}
Assuming Hypothesis~\ref{hyp:weak-GRH}, there is a universal polynomial $h(x,y)$ with the following property.  Let $f \in \Z[x]$ be an irreducible polynomial of degree $d \geq 1$, with leading coefficient $a_d \neq 0$.  Then
\begin{equation} \label{eq:goal-general}
\sum_{\substack{p \leq x \\ p \nmid a_d,\ p\nmid\Delta_g}} W(p) \geq \Li(x) - h(d, \log ||f||) \cdot x^{1-\epsilon}\log x
\end{equation}
for all $x > 2$, where the discriminant $\Delta_g$ of $g(x) = a_d^{d-1}f(x/a_d)$ satisfies \eqref{eq:delta-g}.
\end{proposition}

The polynomial $h(x,y)$ depends on the constant $\epsilon$ in Hypothesis~\ref{hyp:weak-GRH} and the resulting universal constant $c$, but we will not make this explicit other than to note that we can certainly take $h(x,y) = O(x^3y)$.

\begin{proof}[Proof of Theorem~\ref{thm:goal}]
We wish to find $x$, of size polynomial in $d$ and $\log ||f||$, such that the sum on the left side of \eqref{eq:goal-general} is strictly positive.  Proposition~\ref{prop:goal-general} says that this will apply as long as
\[ \Li(x) > h(d,\log ||f||) \cdot x^{1-\epsilon}\log x, \]
and if $x \geq 4$ then the inequality
\[ \Li(x) = \int_2^x \frac{dt}{\log t} \geq \int_2^x \frac{dt}{\log x} = \frac{x - 2}{\log x} \geq \frac{x/2}{\log x} \]
shows that it will further suffice to take $x$ satisfying
\[ x^\epsilon > 2h(d,\log |f||) \cdot (\log x)^2. \]
There is a constant $k > 4$ such that
\[ x^{\epsilon / 3} > \log x \]
for all $x \geq k$, and for all $x \geq \big(2h(d,\log||f||)\big)^{3/\epsilon} + k$ we have
\[ x^\epsilon > x^{\epsilon/3}(\log x)^2 > 2h(d,\log||f||) \cdot (\log x)^2 \]
as desired.
\end{proof}

\bibliographystyle{myalpha}
\bibliography{References}

\end{document}